\documentclass[12pt]{article}
\usepackage[english]{babel}
\usepackage[intlimits]{amsmath}
\usepackage{amsthm,amssymb}
\usepackage{mathrsfs}
\usepackage{indentfirst}
\newcommand{\Beta}{\mathrm{B}}
\DeclareMathOperator{\cov}{cov}
\newtheorem{theorem}{Theorem}

\title{The degree distribution and the number of edges between nodes of given degrees in the Buckley--Osthus model of a random 
web graph}

\author{Evgeniy A. Grechnikov\footnote{Research division in Yandex, Moscow.}}

\date{}

\begin{document}

\maketitle

\renewcommand{\le}{\leqslant}
\renewcommand{\ge}{\geqslant}

\renewcommand{\leq}{\leqslant}
\renewcommand{\geq}{\geqslant}

\begin{abstract}

In this paper, we study some important statistics of the random graph $ H_{a,k}^{(t)} $ in the Buckley--Osthus model, 
where $ t $ is the number of nodes, $ kt $ is the number of edges (so that $ k \in {\mathbb N} $), and $ a > 0 $ 
is the so-called initial attractiveness of a node. This model is a modification of the well-known Bollob\'as--Riordan model. 
First, we find a new asymptotic formula for the expectation of the number $ R(d,t) $ 
of nodes of a given degree $ d $ in a graph in this model. Such a formula is known for $ a \in {\mathbb N} $
and $ d \le t^{1/100(a+1)} $. Both restrictions are unsatisfactory from theoretical and practical points of view. We completely remove 
them. Then we calculate the covariances between any two quantities $ R(d_1,t) $, $ R(d_2,t) $, and using the second moment method we  
show that $ R(d,t) $ is tightly concentrated around its mean for every possible values of $ d $ and $ t $. Furthermore, we study a 
more complicated statistic of the web graph: $ X(d_1,d_2,t) $ is the total number of edges between nodes whose degrees are equal to 
$ d_1 $ and $ d_2 $ respectively. We also find an asymptotic formula for the expectation of $ X(d_1,d_2,t) $ and prove a tight 
concentration result. Again, we do not impose any substantial restrictions on the values of $ d_1, d_2 $, and $ t $. 

\end{abstract}

\section{Introduction}

The real world has many interesting structures which can be thought of as graphs.
A typical example is the World Wide Web: one can consider web pages to be nodes
of a graph and hyperlinks to be edges. One of productive methods for
studying these graphs involves investigation of a suitable random graph model.

First models of random graphs were constructed and investigated
long ago. Classical models and results are systematized, for example,
in \cite{bollobas1} and \cite{Luc}. However, they are not suitable
for approximation of dynamically changing and non-uniform networks.
In particular, the degree sequences of the graphs in these models are very far 
from those observed in reality. 

Recently other models of random graphs were constructed to more closely
match the growth of real networks. One of the first descriptions 
of such a model belongs
to the article \cite{barabasi} by Barab\'asi and Albert.
The authors of this article introduced the ``preferential attachment''
rule. Models following this rule assign the probability of
a new edge to a node according to the current degree of this node, so
more ``popular'' nodes are more attractive for new edges.

However, the article \cite{barabasi} did not contain a
precise model, leaving some parameters unspecified.
Variations of these parameters can significantly change
properties of arising graphs, as shown in \cite{bollobas2},
so one needs something more explicit for theoretical
investigations. Bollob\'as, Riordan et al. proposed an explicit model
in \cite{bollobas3} based on the preferential attachment rule.
In the same article, they rigorously proved a theorem
concerning the degree sequence of a graph in this model. 
Namely, they showed that the number
of nodes with degree $d$ in their model decreases proportional to
$d^{-3}$. The same quantity in real networks decreases
proportional to $d^{-\gamma}$ with different $\gamma$
for different networks, following the so called ``power law''.

The Bollob\'as--Riordan
model has only one parameter, a natural number representing
the ratio of the number of edges to the number of nodes. 
Thus, on the one hand, the Bollob\'as--Riordan model 
does certainly match some real networks
by explaining the power law. But, on the other hand, the number
of parameters in this model is small and does not allow to obtain
the power law with an exponent, which is not equal to $-3$.

In the Bollob\'as--Riordan model, the probability that a node is target for a new edge
is proportional to the degree
of this node. In \cite{dorogovtsev} and \cite{drinea} two groups of researchers independently
proposed to add to the model one more parameter --- an ``initial attractiveness''
of a node which is a positive constant not depending on the degree.
Equivalently, the probability in the proposed model is a linear function
in the degree. However, in the papers \cite{dorogovtsev} and \cite{drinea}, we 
find only some heuristic arguments. 

In \cite{buckley} Buckley and Osthus gave an explicit construction of the above-described model
and rigorously proved a theorem
concerning the degree sequence of a graph in this model when all the parameters
are natural numbers.

Among many articles in this area, we also quote
\cite{jordan}. The model investigated in this article
differs from the Buckley--Osthus model, but the
difference is small, so the results are comparable.
The article deals with the case when parameters
are not necessarily natural.
However, the proven theorem only works for
fixed degree $d$ when the number of nodes tends to infinity;
Bollob\'as et al. as well as Buckley and Osthus allowed
$d$ to grow as some small power of the number of nodes.

There are many other random graph models intended to
approximate real networks. We refer the reader to
\cite{bollobas2} and \cite{durret} for surveys of
such models and corresponding results.

We study the Buckley--Osthus model of a random graph.
Our first goal is to give a significant improvement of the 
above-mentioned theorem
from \cite{buckley} using a completely different method.
We find an asymptotic formula for the expectation of the number
of nodes with degree $d$ without any upper bound
on $d$ and with an estimation of the error term. We also 
prove a tight concentration result.

Since the Bollob\'as--Riordan model is a special case of
the Buckley--Osthus model, our results are also
applicable to it. So, again, we get
a substantial improvement of the main theorem from \cite{bollobas3}.

Our second goal is to study the following quantity.
We fix two numbers $d_1$ and $d_2$. We consider
a node with degree $d_1$ and a node with degree $d_2$.
Then, we calculate the number of edges between these nodes. 
When there are several choices for nodes of given degrees, 
we calculate the mean value.
Since the number of nodes with a fixed degree is known to
have tight concentration around its expectation, it is sufficient to examine
the total number of edges linking a node with degree $d_1$
and a node with degree $d_2$. Here we also obtain an asymptotic 
formula for the expectation and prove a tight concentration result. 

\section{The model and formulation of results}

The Buckley--Osthus model has two parameters, a natural number $k$
and a positive real number $a$.
The number $k$ is the ratio of the number of edges to the number of nodes.
We assume that $a$ and $k$ are constants, so by default all other constants
may depend on them. The Bollob\'as--Riordan model is a special case of
this model with $a=1$.

The model is defined in two stages.
At the first stage, a probability space $H_{a,1}^{(t)}$
is constructed. The elements of $H_{a,1}^{(t)}$ are
undirected graphs with nodes represented by numbers
$1,\dots,t$ and with $t$ edges. The space $H_{a,1}^{(1)}$ contains only
one graph with one node and one loop. Any graph in 
$H_{a,1}^{(t)}$ is obtained from a graph in $H_{a,1}^{(t-1)}$
by adding a new node $t$ and a new edge between
$t$ and a node $\gamma \in \{1, \dots, t\} $ so that
$$
\Pr(\gamma=s)=
\begin{cases}
\frac{\deg_{t-1}(s)-1+a}{(a+1)t-1},&1\le s\le t-1,\\
\frac{a}{(a+1)t-1},&s=t,
\end{cases}
$$
where $\deg_{t-1}$ denotes the degree of a node in the graph from $H_{a,1}^{(t-1)}$.
At the second stage, a final probability space $H_{a,k}^{(t)}$ is
constructed from $H_{a,1}^{(tk)}$ as follows. We take any graph from $H_{a,1}^{(tk)}$. It has 
$ kt $ nodes and $ kt $ edges. We identify the nodes $ 1, \dots, k $; $ k+1, \dots, 2k $; $ \dots $
obtaining $ t $ new nodes, and we keep all the edges obtaining multiple edges and even multiple loops.

We study the number of nodes of degree $d$ in $H_{a,k}^{(t)}$ as a function
of $d$ and $t$. We denote this random quantity by $R(d,t)$
and the value of its expectation by $r(d,t)=ER(d,t)$.

If $d<k$, then clearly $R(d,t)=0$, so it suffices to study the case $d\ge k$. We start by considering $ r(d,t) $. 

\begin{theorem} Let $d\ge k$. The expected value of $R(d,t)$ is
$$
r(d,t)=\frac{\Beta(d-k+ka,a+2)}{\Beta(ka,a+1)}t+O_{a,k}\left(\frac1d\right).
$$
The asymptotic behaviour of the coefficient when $d$ grows is
$$
\frac{\Beta(d-k+ka,a+2)}{\Beta(ka,a+1)}\sim\frac{\Gamma(a+2)}{\Beta(ka,a+1)}d^{-2-a}=(a+1)\frac{\Gamma(ka+a+1)}{\Gamma(ka)}d^{-2-a}.
$$
\end{theorem}

A similar result was obtained in \cite{buckley} (with some factorials instead of Gamma- and Beta-functions). However, for that result, 
it was essential that $ a \in {\mathbb N} $ and $ d \le t^{1/100(a+1)} $. Another result, which can be compared with the one of Theorem 1, 
is proved in \cite{jordan}. It concerns a bit different model, but, nevertheless, it is rather close to our investigations. In this result, 
$ a $ can be any positive real (and so analogous Gamma- and Beta-functions appear in its statement). However, its proof essentially uses 
the assumption that $ d $ is just a constant. In our Theorem 1, we do not have any restrictions on $ d $ and $ a $, and we use a 
completely different method to prove it. 

In fact, Theorem 1 gives an entire picture of what happens to the quantity $ r(d,t) $. If $d=o\left(t^{\frac1{a+1}}\right)$, then Theorem 1
yields the main term of $r(d,t)$. If $d=\Omega\left(t^{\frac1{a+1}}\right)$, then $ r(d,t) $ tends to zero as $ t \to \infty $, which means 
that with high probability there are no nodes of degree $ d $ in a graph in the model. 

Now we want to study in detail the quantity $ R(d,t) $. 

\begin{theorem} Let $d_1\ge k$ and $d_2\ge k$. The covariance between $R(d_1,t)$ and $R(d_2,t)$ is
$$
\cov(R(d_1,t),R(d_2,t))=O_{a,k}\left(\left(d_1^{-2-a}+d_2^{-2-a}\right)t+d_1^{-1}d_2^{-1}\right).
$$
\end{theorem}

Substituting $d_1=d_2=d$ in Theorem 2 and using Chebyshev's inequality, we obtain
the following result.

\vskip+0.2cm

\noindent {\bf Corollary 1.} {\it If $d=d(t)\ge k$ and $\psi(t)\to\infty$ as $t\to\infty$, then
\begin{equation}\label{rconc}
\left|R(d,t)-\frac{\Beta(d-k+ka,a+2)}{\Beta(ka,a+1)}t\right|\le
\left(\sqrt{d^{-a-2}t}+d^{-1}\right)\psi(t)
\end{equation}
with probability tending to 1 as $t\to\infty$.}

\vskip+0.2cm

Let us discuss the meaning of Corollary 1. 

When $d \sim Ct^{\frac1{a+2}}$ with some constant $C$,
both $r(d,t)$ and $\sqrt{d^{-a-2}t}+d^{-1}$ are $O(1)$.
For smaller values of $d$ (i.e., $ d=o\left(t^{\frac1{a+2}}\right) $), inequality (\ref{rconc}) implies the
equivalence (with probability tending to 1 as $d,t\to\infty$)
$$
R(d,t)\sim\frac{(a+1)\Gamma(ka+a+1)}{\Gamma(ka)}d^{-2-a}t.
$$
For larger values of $d$ (i.e., $ t^{\frac1{a+2}}=o(d) $), inequality (\ref{rconc}) means that
$R(d,t)=o(1)$. Since $R(d,t)$ is an integer number by definition, $ R(d,t)=0 $
(again, with probability tending to 1 as $d,t\to\infty$). Thus, we have an almost entire picture of what happens to $ R(d,t) $. 

\vskip+0.2cm

We also study the total number of edges linking a node with
degree $d_1$ and a node with degree $d_2$. We denote
this random quantity by $X(d_1,d_2,t)$. When $d_1=d_2$,
we count every edge twice, but do not count loops.

\begin{theorem} Let $d_1\ge k$ and $d_2\ge k$. There exists
a function $c_X(d_1,d_2)$ such that
$$
EX(d_1,d_2,t)=c_X(d_1,d_2)t+O_{a,k}(1)
$$
and
\begin{multline*}
c_X(d_1,d_2)=
\frac{\Gamma(d_1-k+ka)\Gamma(d_2-k+ka)\Gamma(d_1+d_2-2k+2ka+3)}{\Gamma(d_1-k+ka+2)\Gamma(d_2-k+ka+2)\Gamma(d_1+d_2-2k+2ka+a+2)} \times \\ \times
ka(a+1)\frac{\Gamma(ka+a+1)}{\Gamma(ka)}\left(1+\theta(d_1,d_2)\frac{(d_1-k+ka+1)(d_2-k+ka+1)}{(d_1+d_2-2k+2ka+1)(d_1+d_2-2k+2ka+2)}\right),
\end{multline*}
where
$$
-4+\frac2{1+ka}\le\theta(d_1,d_2)\le a\frac{\Gamma(ka+1)\Gamma(2ka+a+3)}{\Gamma(2ka+2)\Gamma(ka+a+2)}.
$$
When both $d_1$ and $d_2$ grow, the asymptotic behaviour of $c_X$ is
$$
c_X(d_1,d_2)=ka(a+1)\frac{\Gamma(ka+a+1)}{\Gamma(ka)}\frac{(d_1+d_2)^{1-a}}{d_1^2d_2^2}\left(1+O_{a,k}\left(\frac1{d_1}+\frac1{d_2}+\frac{d_1d_2}{(d_1+d_2)^2}\right)\right).
$$
\end{theorem}
Note that the last formula in Theorem 3 does not give an asymptotic behaviour if $d_1$ and $d_2$
grow so that $\frac{d_2}{d_1}$ tends to a finite nonzero limit. The precise bounds show that
the term $\frac{(d_1+d_2)^{1-a}}{d_1^2d_2^2}$ still gives the correct order of growth for $c_X$,
but the coefficient can differ from $ka(a+1)\frac{\Gamma(ka+a+1)}{\Gamma(ka)}$. And in fact, the
coefficient differs.

\begin{theorem} Let $d_1,d_2\ge k$, $c>0$. Then
$$
P\left(|X(d_1,d_2,t)-EX(d_1,d_2,t)|\ge c(d_1+d_2)\sqrt{kt}\right)\le2\exp\left(-\frac{c^2}8\right).
$$
In particular, if $c(t)\to\infty$ as $t\to\infty$, then $|X-EX|<c(t)(d_1+d_2)\sqrt{kt}$
with probability tending to 1.
\end{theorem}

From Theorems 3 and 4, we immediately obtain the following assertion.

\vskip+0.2cm

\noindent {\bf Corollary 2.} {\it If $ (d_1+d_2)^ad_1^2d_2^2=o\left(\sqrt{t}\right) $, then 
with probability tending to 1 as $d_1,d_2,t\to\infty$
$$
X(d_1,d_2,t) \sim c_X(d_1,d_2)t. 
$$}

The mean value of the number of edges between one node with degree $d_1$ and
another node with degree $d_2$ is $\frac{X(d_1,d_2,t)}{R(d_1,t)R(d_2,t)}$.
Since the quantities $R(d,t)$ and $X(d_1,d_2,t)$ are tightly concentrated around their expectations,
the main term of the ratio is
$$
\frac{\Gamma(ka+1)}{(a+1)\Gamma(ka+a+1)}\frac{d_1^ad_2^a(d_1+d_2)^{1-a}}t.
$$
Again, the constant factor can differ if $d_1$ and $d_2$ grow so that $\frac{d_2}{d_1}$ tends
to a finite nonzero limit, but the order is correct even in this case.

\section{Proof of Theorem 1}
For a property $P$, we denote
$$
[P]=\begin{cases}1,&P\mbox{ holds},\\0,&\mbox{otherwise}.\end{cases}
$$

First of all, we reformulate the model without references to $H_{a,1}^{(t)}$. The probability
space $H_{a,k}^{(1)}$ obviously consists of one graph with 1 node and $k$ loops.
The space $H_{a,k}^{(t+1)}$ can be obtained from $H_{a,k}^{(t)}$ by adding to any graph from 
$H_{a,k}^{(t)}$ a new node $t+1$ and $k$ edges in the following $ k $ steps.
At the $i$th step, we add one edge between the new
node and one of the existing nodes $\gamma$. If $\gamma\ne t+1$,
then it corresponds to a group of nodes $\gamma_1,\dots,\gamma_k$ in $H_{a,1}^{(kt+i-1)}$.
The sum of degrees of $\gamma_1,\dots,\gamma_k$ equals the degree of $\gamma$ in the
graph before the $i$th step. We denote this degree by $\deg_{t,i}$. So
$$\Pr(\gamma=s)=\frac{\deg_{t,i}(s)+k(a-1)}{(a+1)(kt+i)-1},\quad 1\le s\le t.$$

If $\gamma=t+1$, the corresponding group in $H_{a,1}^{(kt+i-1)}$
has only $i-1$ nodes. Hence,
$$\Pr(\gamma=t+1)=\frac{\deg_{t,i}(t+1)+(i-1)(a-1)+a}{(a+1)(kt+i)-1}.$$

We want to express any value $r(d,t)$ in terms of some values with smaller $t$.
Let us consider the transition from $H_{a,k}^{(t)}$ to
$H_{a,k}^{(t+1)}$. Let $r(d,t,i)$ denote the average number of nodes of degree $d$,
not including the last node $t+1$, before the $i$th step, and $r(d,t,i+1)$ --- the similar
number after the $i$th step. Let $\gamma$ be a head of the edge added in the $i$th step. Then, 
\begin{multline}\label{rirecur}
r(d,t,i+1)=\sum_{s=1}^{t}\Pr(\deg_{t,i+1}(s)=d)=
\sum_{s=1}^{t}(\Pr(\deg_{t,i+1}(s)=d,\gamma=s)+\\+
\Pr(\deg_{t,i+1}(s)=d,\gamma\ne s))=
\sum_{s=1}^{t}(\Pr(\deg_{t,i}(s)=d-1,\gamma=s)+\\+
\Pr(\deg_{t,i}(s)=d,\gamma\ne s))=
\sum_{s=1}^{t}\Bigg(\Pr(\deg_{t,i}(s)=d-1)\frac{(d-1)+k(a-1)}{(a+1)(kt+i)-1}+\\+
\Pr(\deg_{t,i}(s)=d)\left(1-\frac{d+k(a-1)}{(a+1)(kt+i)-1}\right)\Bigg)=\\=
r(d-1,t,i)\frac{(d-1)+k(a-1)}{(a+1)(kt+i)-1}+
r(d,t,i)\left(1-\frac{d+k(a-1)}{(a+1)(kt+i)-1}\right).
\end{multline}
By definition,
\begin{equation}\label{riandr}
r(d,t)=r(d,t,1),\qquad r(d,t+1)=r(d,t,k+1)+\Pr(\deg_{t,k+1}(t+1)=d).
\end{equation}
The function $r(d,t)$ is completely determined by the equations
(\ref{rirecur}), (\ref{riandr}) and the starting condition
\begin{equation}\label{rstart}
r(d,1)=[d=2k].
\end{equation}

The equation (\ref{riandr}) includes the function $\Pr(\deg_{t,k+1}(t+1)=d)$.
Obviously,
\begin{equation}\label{probzero}
\Pr(\deg_{t,k+1}(t+1)=d)=0,\quad d<k\mbox{ or }d>2k.
\end{equation}
The minimal value $\deg_{t,k+1}(t+1)=k$ is obtained when no one of the $k$ edges
is a loop. In this case, $\deg_{t,i}(t+1)=i-1$ for all $i$, so
$$
\Pr(\deg_{t,k+1}(t+1)=k)=\prod_{i=1}^k\left(1-\frac{ia}{(a+1)(kt+i)-1}\right)=1+O\left(\frac1t\right).
$$
(Note that a constant in $O()$ depends on $a$ and $k$).


Because $\sum\limits_{d=k}^{2k}\Pr(\deg_{t,k+1}=d)=1$ and $\Pr(\deg_{t,k+1}=d)\ge 0$, we get
$$
\Pr(\deg_{t,k+1}(t+1)=d)=O\left(\frac1t\right),\quad k<d\le 2k.
$$
Since $d$ is bounded in the last equality, its right hand side can be equivalently
written as $O\left(\frac1{d^2t}\right)$.

Let
\begin{equation}\label{cdef}
c(d)=\begin{cases}\frac{\Beta(d-k+ka,a+2)}{\Beta(ka,a+1)},&d\ge k,\\0,&d<k.\end{cases}
\end{equation}

If $d>k$, then
\begin{multline*}
\frac{c(d-1)}{c(d)}=\frac{\Beta(d-1-k+ka,a+2)}{\Beta(d-k+ka,a+2)}=\\=
\frac{\Gamma(d-1-k+ka)/\Gamma(d+1-k+ka+a)}{\Gamma(d-k+ka)/\Gamma(d-k+ka+a+2)}=
\frac{d+1-k+ka+a}{d-1-k+ka}.
\end{multline*}
Also
$$
c(k)=\frac{\Beta(ka,a+2)}{\Beta(ka,a+1)}=\frac{\Gamma(a+2)/\Gamma(ka+a+2)}{\Gamma(a+1)/\Gamma(ka+a+1)}=
\frac{a+1}{ka+a+1}.
$$
In particular, $c(d-1)>c(d)$, so $c(d)<c(k)<1$ for all $d\ge k$.

For the rest of the proof, we will assume that $d\ge k$. Note that of course this does not
imply $d-1\ge k$.

When $d$ grows, the asymptotic behaviour of $c(d)$ is
\begin{multline*}
\ln c(d)=\ln\frac{\Gamma(a+2)}{\Beta(ka,a+1)}\frac{\Gamma(d-k+ka)}{\Gamma(d-k+ka+a+2)}=
\ln\frac{\Gamma(a+2)}{\Beta(ka,a+1)}+\\+
(d-k+ka)\left(\ln(d-k+ka)-1\right)-(d-k+ka+a+2)\left(\ln(d-k+ka+a+2)-1\right)+O\left(\frac 1d\right)=\\=
\ln\frac{\Gamma(a+2)}{\Beta(ka,a+1)}+(d-k+ka)(\ln d+\frac{-k+ka}d-1)-(d-k+ka+a+2)(\ln d+\frac{-k+ka+a+2}d-1)\\+
O\left(\frac 1d\right)=
\ln\frac{\Gamma(a+2)}{\Beta(ka,a+1)}-(a+2)\ln d+O\left(\frac 1d\right),
\end{multline*}
\begin{equation}\label{casymp}
c(d)=\frac{\Gamma(a+2)}{\Beta(ka,a+1)}d^{-2-a}\left(1+O\left(\frac 1d\right)\right).
\end{equation}

Let
$$
\tilde r(d,t,i)=r(d,t,i)-c(d)\left(t+\frac ik-\frac1{k(a+1)}\right).
$$
It is easy to see that the theorem is equivalent to $\tilde r(d,t,i)=O(1)$.
Using (\ref{rirecur}), we obtain
\begin{multline}\label{ritilderecur}
\tilde r(d-1,t,i)\frac{(d-1)+k(a-1)}{(a+1)(kt+i)-1}+
\tilde r(d,t,i)\left(1-\frac{d+k(a-1)}{(a+1)(kt+i)-1}\right)=\\=
r(d,t,i+1)-
\left(t+\frac ik-\frac1{k(a+1)}\right)
\Bigg(c(d-1)\frac{(d-1)+k(a-1)}{(a+1)(kt+i)-1}+\\+
c(d)\left(1-\frac{d+k(a-1)}{(a+1)(kt+i)-1)}\right)\Bigg)=
r(d,t,i+1)-
\left(t+\frac ik-\frac1{k(a+1)}\right)c(d)\times\\ \times
\frac{(1-[d=k])(d+1-k+ka+a)+(a+1)(kt+i)-1-(d+k(a-1))}{(a+1)(kt+i)-1}=\\=
r(d,t,i+1)-
c(d)\frac{(a+1)(kt+i+1)-1-[d=k](1+ka+a)}{k(a+1)}=\\=
\tilde r(d,t,i+1)+[d=k]c(k)\frac{1+ka+a}{k(a+1)}=\tilde r(d,t,i+1)+\frac{[d=k]}k.
\end{multline}

Let $C=C(a,k)$ be a sufficiently large constant which will be determined later. We claim
that
\begin{equation}\label{ritildebound}
\left|\tilde r(d,t,i)+(i-1)\frac{[d=k]}k\right|\le\frac{C}{d+ka}\left(1-\frac{\min\{1,ka\}}{(a+2)(t+1)(d+ka)}\right)^{i-1}
\end{equation}
for all $i=1,\dots,k+1$ and for all natural $d\ge k$ and $t$. Note that this implies
$\tilde r(d,t,i)=O\left(\frac1d\right)$ and Theorem 1.

The equations (\ref{rirecur}), (\ref{riandr}), (\ref{rstart}), (\ref{probzero}) imply that
$r(d,t,i)=0$ if $d>kt+i-1+k$. In this case, using (\ref{casymp}), we obtain
$$
\tilde r(d,t,i)=-c(d)\left(t+\frac ik-\frac1{k(a+1)}\right)=O\left(d^{-2-a}t\right)=O\left(d^{-1-a}\right),
$$
so if $d>kt+i-1+k$, (\ref{ritildebound}) is true for all sufficiently large values of $C$.

Now assume $d\le kt+i-1+k$. We will prove (\ref{ritildebound}) by induction on $t$ and, for fixed $t$, on $i$.
The basis of induction $t=1,\dots,1+\lfloor\frac1{ka}\rfloor$
and any $i=1,\dots,k+1$ obviously holds for all sufficiently large values of $C$.

Now let $t\ge2+\lfloor\frac1{ka}\rfloor$ and let (\ref{ritildebound}) hold for $t-1$. Using (\ref{riandr}), we obtain
$$
\tilde r(d,t,1)=\tilde r(d,t-1,k+1)+\Pr(\deg_{t-1,k+1}(t)=d)=
\tilde r(d,t-1,k+1)+[d=k]+O\left(\frac 1{d^2t}\right).
$$
Therefore, 
\begin{multline*}
|\tilde r(d,t,1)|\le\frac{C}{d+ka}\left(1-\frac{\min\{1,ka\}}{(a+2)(t+1)(d+ka)}\right)^k+O\left(\frac1{d^2t}\right) \le \\ \le
\frac{C}{d+ka}-\frac{C\min\{1,ka\}/(a+2)}{(t+1)(d+ka)^2}+O\left(\frac1{d^2t}\right).
\end{multline*}
Thus, induction step on $t$ is proved for all sufficiently large values of $C$.

Finally, let $t\ge2+\lfloor\frac1{ka}\rfloor>1+\frac1{ka}$, $i>1$ and let (\ref{ritildebound}) hold for $i-1$.
We temporarily denote $T=(a+1)(kt+i-1)-1$. Note that $T$ depends on $t$ and $i$, but not on $d$.
From (\ref{ritilderecur}) we obtain
\begin{multline*}
\tilde r(d,t,i)+(i-1)\frac{[d=k]}k=(i-2)\frac{[d=k]}k+\\+
\tilde r(d-1,t,i-1)\frac{(d-1)+k(a-1)}{T}+
\tilde r(d,t,i-1)\left(1-\frac{d+k(a-1)}{T}\right)=\\=
\left(\tilde r(d-1,t,i-1)+(i-2)\frac{[d-1=k]}k\right)\frac{(d-1)+k(a-1)}{T}+\\+
\left(\tilde r(d,t,i-1)+(i-2)\frac{[d=k]}k\right)\left(1-\frac{d+k(a-1)}{T}\right)+[k\le d\le k+1]O\left(\frac1t\right).
\end{multline*}
The remainder term $[k\le d\le k+1]O\left(\frac1t\right)$ can be written as $O\left(\frac1{d^2t}\right)$.
The assumptions $d\le kt+i-1+k$ and $t>1+\frac1{ka}$ imply that $1-\frac{d+k(a-1)}{(a+1)(kt+i-1)-1}\ge0$.
If $d>k$, then
\begin{multline*}
\frac{C}{d-1+ka}\frac{(d-1)+k(a-1)}{T}+\frac{C}{d+ka}\left(1-\frac{d+k(a-1)}{T}\right)=\\=
\frac{C}{d+ka}\left(1-\frac{1}{T}\left(-\frac{d+ka}{d-1+ka}((d-1)+k(a-1))+d+k(a-1)\right)\right)=\\=
\frac{C}{d+ka}\left(1-\frac{1}{T}\left(1-\frac{(d-1)+k(a-1)}{d-1+ka}\right)\right)=
\frac{C}{d+ka}\left(1-\frac{k}{T(d-1+ka)}\right) \le \\ \le
\frac{C}{d+ka}\left(1-\frac{k\min\{1,ka\}}{T(d+ka)}\right).
\end{multline*}
If $d=k$, then
$$
\frac{C}{d+ka}\left(1-\frac{d+k(a-1)}T\right)=
\frac{C}{d+ka}\left(1-\frac{ka(k+ka)}{T(d+ka)}\right)\le\frac{C}{d+ka}\left(1-\frac{k\min\{1,ka\}}{T(d+ka)}\right).
$$
In both cases,
\begin{multline}\label{ritildetmp}
\left|\tilde r(d,t,i)+(i-1)\frac{[d=k]}k\right|\le
\frac{C}{d+ka}\left(1-\frac{\min\{1,ka\}}{(a+2)(t+1)(d+ka)}\right)^{i-2}\left(1-\frac{k\min\{1,ka\}}{T(d+ka)}\right)+\\+
O\left(\frac1{d^2t}\right).
\end{multline}
Note that
\begin{multline*}
1-\frac{k\min\{1,ka\}}{((a+1)(kt+i-1)-1)(d+ka)}=\left(1-\frac{k\min\{1,ka\}}{(a+2)(kt+k)(d+ka)}\right)-\\-
\frac{k\min\{1,ka\}}{(a+1)(a+2)kt(d+ka)}\left(1+O\left(\frac1t\right)\right),
\end{multline*}
where the left part is strictly less than the first term in the right part, so that
the second term in the right part is positive. Now (\ref{ritildetmp}) implies (\ref{ritildebound})
for all sufficiently large values of $C$.

Theorem 1 is proved. 

\section{Proof of Theorem 2}
By definition and linearity of expectation,
\begin{multline}\label{covmain}
\cov(R(d_1,t),R(d_2,t))=E(R(d_1,t)R(d_2,t))-r(d_1,t)r(d_2,t)=\\=
E\sum_{s_1,s_2=1}^t[\deg s_1=d_1,\deg s_2=d_2]-r(d_1,t)r(d_2,t)=\\=
\sum_{s_1\ne s_2}\Pr(\deg s_1=d_1,\deg s_2=d_2)+[d_1=d_2]r(d_1,t)-r(d_1,t)r(d_2,t).
\end{multline}
We will estimate the sum
$$
r_2(d_1,d_2,t)=\sum_{\substack{s_1,s_2=1\\s_1\ne s_2}}^t \Pr(\deg s_1=d_1,\deg s_2=d_2)
$$
as we have done it for the function $r(d,t)$ in the proof of Theorem 1.

As with $r(d,t)$, we define a function $r_2(d_1,d_2,t,i)$ as the value of $r_2(d_1,d_2,t)$ before the $i$th step
in the transition from $H_{a,k}^{(t)}$ to $H_{a,k}^{(t+1)}$. The recurrent equation is deduced similarly to
(\ref{rirecur}). For fixed $s_1$ and $s_2$,
there are three non-intersecting cases: $\gamma=s_1$, $\gamma=s_2$, and $\gamma\not\in\{s_1,s_2\}$.
In the first case, we get
\begin{multline*}
\Pr(\deg_{t,i+1}(s_1)=d_1,\deg_{t,i+1}(s_2)=d_2,\gamma=s_1)=\\=
\Pr(\deg_{t,i}(s_1)=d_1-1,\deg_{t,i+1}(s_2)=d_2,\gamma=s_1)=\\=
\Pr(\deg_{t,i}(s_1)=d_1-1,\deg_{t,i+1}(s_2)=d_2)\frac{d_1-1+k(a-1)}{(a+1)(kt+i)-1}.
\end{multline*}
The second case is the same with $d_1$ and $d_2$ interchanged. In the third case, we get
\begin{multline*}
\Pr(\deg_{t,i+1}(s_1)=d_1,\deg_{t,i+1}(s_2)=d_2,\gamma\ne s_1,\gamma\ne s_2)=\\=
\Pr(\deg_{t,i}(s_1)=d_1,\deg_{t,i+1}(s_2)=d_2,\gamma\ne s_1,\gamma\ne s_2)=\\=
\Pr(\deg_{t,i}(s_1)=d_1,\deg_{t,i+1}(s_2)=d_2)\left(1-\frac{d_1+k(a-1)+d_2+k(a-1)}{(a+1)(kt+i)-1}\right),
\end{multline*}
so the final formula is
\begin{multline}\label{r2recur}
r_2(d_1,d_2,t,i+1)=r_2(d_1-1,d_2,t,i)\frac{(d_1-1)+k(a-1)}{(a+1)(kt+i)-1}+\\+
r_2(d_1,d_2-1,t,i)\frac{(d_2-1)+k(a-1)}{(a+1)(kt+i)-1}+\\+
r_2(d_1,d_2,t,i)\left(1-\frac{d_1+d_2+2k(a-1)}{(a+1)(kt+i)-1}\right). 
\end{multline}
By definition,
$$
r_2(d_1,d_2,t)=r_2(d_1,d_2,t,1),
$$
\begin{multline}\label{r2andri}
r_2(d_1,d_2,t+1)=r_2(d_1,d_2,t,k+1)+\\+
\sum_{s=1}^{t}\Pr(\deg_{t,k+1}(s)=d_1,\deg_{t,k+1}(t+1)=d_2)+\\+
\sum_{s=1}^{t}\Pr(\deg_{t,k+1}(s)=d_2,\deg_{t,k+1}(t+1)=d_1),
\end{multline}
and the starting condition is
$$r_2(d_1,d_2,1)=0.$$

The equation (\ref{r2recur}) includes a function
\begin{equation}\label{r2primedef}
r_2'(d_1,d_2,t,i)=\sum_{s=1}^{t}\Pr(\deg_{t,i}(s)=d_1,\deg_{t,i}(t)=d_2)
\end{equation}
(and the same function with swapped arguments), so we will first estimate $r_2'$.
Again, we write a recurrent equation. For fixed $s$, there are three non-intersecting cases:
$\gamma=s$, $\gamma=t+1$, $\gamma\not\in\{s,t+1\}$. If $\deg_{t,i+1}(s)=d_1$, then
$\deg_{t,i}(s)$ equals $d_1-1$ in the first case and $d_1$ in the two other cases.
If $\deg_{t,i+1}(t+1)=d_2$, then $\deg_{t,i}(t+1)=d_2-2$ in the second case
and $d_2-1$ in the two other cases. Calculating probabilities, we obtain
\begin{multline*}
r_2'(d_1,d_2,t,i+1)=r_2'(d_1-1,d_2-1,t,i)\frac{(d_1-1)+k(a-1)}{(a+1)(kt+i)-1}+\\+
r_2'(d_1,d_2-2,t,i)\frac{(d_2-2)+(i-1)(a-1)+a}{(a+1)(kt+i)-1}+\\+
r_2'(d_1,d_2-1,t,i)\left(1-\frac{d_1+k(a-1)+(d_2-1)+(i-1)(a-1)+a}{(a+1)(kt+i)-1}\right).
\end{multline*}

Before the $1$st step, the node $t$ has degree 0, so
$$
r_2'(d_1,d_2,t,1)=\sum_{s=1}^{t}\Pr(\deg_{t,1}(s)=d_1)[d_2=0]=[d_2=0]r(d_1,t).
$$
We continue to use notation (\ref{cdef}). Obviously, $r_2'(d_1,d_2,t,i)=0$ when
$d_2 > 2(i-1)$ or $d_1\ge2(kt+i)$. If $d_1<2(kt+i)$ and $d_2\le2(i-1)$,
then $\frac{d_1+d_2}t\cdot O(d_1^{-a-2}t)=O(d_1^{-a-1})=O(d_1^{-1})$ and $\frac{d_1+d_2}t\cdot O(d_1^{-1})=O(d_1^{-1})$.
Now it is easy to see that
\begin{equation}\label{r2prime}
r_2'(d_1,d_2,t,i)=[d_2=i-1]c(d_1)t+O(d_1^{-1}).
\end{equation}
Since $r_2'(d_1,d_2,t,i)=0$ when $d_2\ge 2k$, the remainder term in (\ref{r2prime})
is zero when $d_2\ge 2k$ and can be written as $O\left(\frac{1}{d_1d_2^{a+1}}\right)$.

Since $r(d_1,t)=0$ when $d_1>kt+k$, $r_2'(d_1,d_2,t,i)=0$
when $d_1+d_2>kt+k+2(i-1)$. In particular,
$r_2'(d_1,d_2,t,k+1)+r_2'(d_2,d_1,t,k+1)=0$ when
$d_1+d_2>kt+2k$. By (\ref{r2recur}) and (\ref{r2andri}),
$r_2(d_1,d_2,t,i)=0$ when $d_1+d_2>kt+2k+(i-1)$.

For the rest of the proof, we will assume that $d_1\ge k$ and $d_2\ge k$.
Note that of course this does not imply $d_1-1\ge k$ and $d_2-1\ge k$.

Let
$$
\tilde r_2(d_1,d_2,t,i)=r_2(d_1,d_2,t,i)-
c(d_1)c(d_2)
\left(t+\frac ik-\frac1{k(a+1)}\right)\left(t+\frac{i+1}k-\frac1{k(a+1)}\right).
$$
We temporarily denote $T=kt+i-\frac1{a+1}$.
We express $r_2$ in terms of $\tilde r_2$ and use
(\ref{r2recur}). In the expression, there are terms with $\tilde r_2$ with various arguments.
Now we transform the terms without $\tilde r_2$ from the right part
of (\ref{r2recur}).
\begin{multline*}
\frac Tk\frac{T+1}k\Bigg(
c(d_1-1)c(d_2)\frac{(d_1-1)+k(a-1)}{(a+1)T}+\\+
c(d_1)c(d_2-1)\frac{(d_2-1)+k(a-1)}{(a+1)T}+
c(d_1)c(d_2)\left(1-\frac{d_1+d_2+2k(a-1)}{(a+1)T}\right)
\Bigg)=\\=
\frac Tk\frac{T+1}kc(d_1)c(d_2)\Bigg(
(1-[d_1=k])\frac{d_1+1-k+ka+a}{(a+1)T}+\\+
(1-[d_2=k])\frac{d_2+1-k+ka+a}{(a+1)T}+
1-\frac{d_1+d_2+2k(a-1)}{(a+1)T}\Bigg)=\\=
\frac Tk\frac{T+1}kc(d_1)c(d_2)\Bigg(
\frac{2(a+1)+(a+1)T-[d_1=k](1+ka+a)-[d_2=k](1+ka+a)}{(a+1)T}
\Bigg)=\\=
\frac{T+1}k\frac{T+2}kc(d_1)c(d_2)-[d_1=k]\frac{(T+1)c(d_2)}{k^2}-[d_2=k]\frac{(T+1)c(d_1)}{k^2}.
\end{multline*}
The first term equals the term without $\tilde r_2$ in the left part, so
\begin{multline*}
\tilde r_2(d_1,d_2,t,i+1)=\tilde r_2(d_1-1,d_2,t,i)\frac{(d_1-1)+k(a-1)}{(a+1)(kt+i)-1}+\\+
\tilde r_2(d_1,d_2-1,t,i)\frac{(d_2-1)+k(a-1)}{(a+1)(kt+i)-1}+
\tilde r_2(d_1,d_2,t,i)\left(1-\frac{d_1+d_2+2k(a-1)}{(a+1)(kt+i)-1}\right)-\\-
[d_1=k]\frac{\left(kt+i-\frac1{a+1}+1\right)c(d_2)}{k^2}-[d_2=k]\frac{\left(kt+i-\frac1{a+1}+1\right)c(d_1)}{k^2}.
\end{multline*}
Relations (\ref{r2andri}) and (\ref{r2prime}) imply
$$
\tilde r_2(d_1,d_2,t,1)=\tilde r_2(d_1,d_2,t-1,k+1)+[d_2=k]c(d_1)t+[d_1=k]c(d_2)t+
O\left(\frac{1}{d_1^{a+1}d_2}+\frac{1}{d_1d_2^{a+1}}\right).
$$

Let
$$
c_1(d_1',d_2')=\begin{cases}\displaystyle\frac{\Gamma(d_1'-k+ka)}{(d_2'+k(a-1))\Gamma(d_1'-k+ka+a+1)},&d_1'\ge k,d_2'\ge k,\\
0,&d_1'<k\mbox{ or }d_2'<k.\end{cases}
$$
By definition, for $d_1>k$,
$$
\frac{c_1(d_1-1,d_2)}{c_1(d_1,d_2)}=\frac{d_1-k+ka+a}{d_1-k+ka-1}.
$$
For $d_2>k$,
$$
\frac{c_1(d_1,d_2-1)}{c_1(d_1,d_2)}=\frac{d_2+k(a-1)}{d_2-1+k(a-1)}.
$$
Similar to (\ref{casymp}),
$$
c_1(d_1,d_2)=\frac{d_1^{-1-a}}{d_2+k(a-1)}\left(1+O\left(\frac1{d_1}\right)\right).
$$
Moreover,
\begin{multline*}
\frac{c_1(d_1,d_2)}{c(d_1)}=
\frac{\Beta(ka,a+1)}{\Gamma(a+2)}\frac{\Gamma(d_1-k+ka)\Gamma(d_1-k+ka+a+2)}{(d_2+k(a-1))\Gamma(d_1-k+ka+a+1)\Gamma(d_1-k+ka)}=\\=
\frac{\Beta(ka,a+1)}{\Gamma(a+2)}\frac{d_1-k+ka+a+1}{d_2+k(a-1)}.
\end{multline*}

Let $C=C(a,k)$ be a sufficiently large constant which will be determined later. We claim that
\begin{multline}\label{r2tildebound}
\left|\tilde r_2(d_1,d_2,t,i)+\frac{i-1}k([d_2=k]c(d_1)+[d_1=k]c(d_2))t\right|\le\\ \le
C(c_1(d_1,d_2)+c_1(d_2,d_1))\left(kt+\frac{a+\frac12}{a+1}i\right)
\end{multline}
for all $i=1,\dots,k+1$ and for all natural $d_1\ge k,d_2\ge k,t$. Since both parts of (\ref{r2tildebound})
are symmetric in $d_1$ and $d_2$, it is sufficient to consider the case $d_1\le d_2$.

If $d_1+d_2>kt+2k+(i-1)$, then $r_2(d_1,d_2,t,i)=0$ and
\begin{multline*}
\frac{\left|\tilde r_2(d_1,d_2,t,i)\right|}{c_1(d_1,d_2)t}=
\frac{c(d_1)c(d_2)\left(t+\frac ik-\frac1{k(a+1)}\right)\left(t+\frac{i+1}k-\frac1{k(a+1)}\right)}{c_1(d_1,d_2)t}=\\=
O\left(\frac{d_2+k(a-1)}{d_1-k+ka+a+1}c(d_2)t\right)=O\left(\frac{t}{d_1d_2^{1+a}}\right).
\end{multline*}
Since $d_2\ge\frac{d_1+d_2}2\ge\frac k2t+k$, the right part is bounded. Obviously, $[d_2=k]=0$
and $[d_1=k]c(d_2)t=O\left(\frac{t}{d_2^{2+a}}\right)$ is bounded too. Thus
(\ref{r2tildebound}) holds when $d_1+d_2>kt+2k+(i-1)$ for all sufficiently large values of $C$.

Let $d_1+d_2\le kt+2k+(i-1)$. We will prove (\ref{r2tildebound}) by induction on $t$ and, for fixed $t$,
on $i$. The basis of induction $t=1,\dots,2+\lfloor\frac1{ka}\rfloor$ and any $i=1,\dots,k+1$ obviously holds for all sufficiently large
values of $C$.

Let $t\ge3+\lfloor\frac1{ka}\rfloor$ and let (\ref{r2tildebound}) hold for $t-1$. We continue to use the restriction $d_1\le d_2$. Thus, 
\begin{multline*}
|\tilde r_2(d_1,d_2,t,1)|=
\left|\tilde r_2(d_1,d_2,t-1,k+1)+[d_2=k]c(d_1)t+[d_1=k]c(d_2)t+O\left(\frac{1}{d_1^{a+1}d_2}\right)\right|\le\\ \le
C(c_1(d_1,d_2)+c_1(d_2,d_1))\left(k(t-1)+\frac{a+\frac12}{a+1}(k+1)\right)+O\left(\frac{1}{d_1^{a+1}d_2}\right).
\end{multline*}
Since $c_1(d_1,d_2)=O\left(\frac1{d_1^{a+1}d_2}\right)$, the right part is less than
$C(c_1(d_1,d_2)+c_1(d_2,d_1))\left(kt+\frac{a+\frac12}{a+1}\right)$
for all sufficiently large values of $C$. This completes the induction on $t$.

Let $t\ge3+\lfloor\frac1{ka}\rfloor>2+\frac1{ka}$, $i>1$ and let (\ref{r2tildebound}) hold for $i-1$. Then, 
$$
\tilde r_2(d_1,d_2,t,i)+\frac{i-1}k([d_2=k]c(d_1)+[d_1=k]c(d_2))t=
$$
$$
=\tilde r_2(d_1-1,d_2,t,i-1)\frac{(d_1-1)+k(a-1)}{(a+1)(kt+i-1)-1}+
\tilde r_2(d_1,d_2-1,t,i-1)\frac{(d_2-1)+k(a-1)}{(a+1)(kt+i-1)-1}+
$$
$$
+\tilde r_2(d_1,d_2,t,i-1)\left(1-\frac{d_1+d_2+2k(a-1)}{(a+1)(kt+i-1)-1}\right)-
[d_1=k]\frac{\left(kt+i-1-\frac1{a+1}+1\right)c(d_2)}{k^2}-
$$
$$
-[d_2=k]\frac{\left(kt+i-1-\frac1{a+1}+1\right)c(d_1)}{k^2}+
\frac{i-1}k([d_2=k]c(d_1)+[d_1=k]c(d_2))t=
$$
\begin{multline*}
=\left(\tilde r_2(d_1\!-\!1,d_2,t,i-1)+\frac{i-2}k([d_2=k]c(d_1\!-\!1)+[d_1\!-\!1=k]c(d_2))t\right)\frac{(d_1\!-\!1)+k(a-1)}{(a+1)(kt+i-1)-1}+\\+
\left(\tilde r_2(d_1,d_2\!-\!1,t,i-1)+\frac{i-2}k([d_2\!-\!1=k]c(d_1)+[d_1=k]c(d_2\!-\!1))t\right)\frac{(d_2\!-\!1)+k(a-1)}{(a+1)(kt+i-1)-1}+\\+
\left(\tilde r_2(d_1,d_2,t,i-1)+\frac{i-2}k([d_2=k]c(d_1)+[d_1=k]c(d_2))t\right)\left(1-\frac{d_1+d_2+2k(a-1)}{(a+1)(kt+i-1)-1}\right)+\\+
[d_1\le k+1]O(c(d_2))+[d_2\le k+1]O(c(d_1)).
\end{multline*}
The assumptions $d_1+d_2\le kt+2k+(i-1)$ and $t>2+\frac1{ka}$ imply that $1-\frac{d_1+d_2+2k(a-1)}{(a+1)(kt+i-1)-1}\ge0$. By the induction 
hypothesis
\begin{multline*}
\left|\tilde r_2(d_1,d_2,t,i)+\frac{i-1}k([d_2=k]c(d_1)+[d_1=k]c(d_2))t\right|\le \\ \le
C\left(kt+\frac{a+\frac12}{a+1}(i-1)\right)\Bigg((c_1(d_1-1,d_2)+c_1(d_2,d_1-1))\frac{(d_1-1)+k(a-1)}{(a+1)(kt+i-1)-1}+\\+
(c_1(d_1,d_2-1)+c_1(d_2-1,d_1))\frac{(d_2-1)+k(a-1)}{(a+1)(kt+i-1)-1}+\\+
(c_1(d_1,d_2)+c_1(d_2,d_1))\left(1-\frac{d_1+d_2+2k(a-1)}{(a+1)(kt+i-1)-1}\right)\Bigg)+\\+
[d_1\le k+1]O(c(d_2))+[d_2\le k+1]O(c(d_1)).
\end{multline*}
Since
\begin{multline*}
\frac{c_1(d_1-1,d_2)}{c_1(d_1,d_2)}\frac{(d_1-1)+k(a-1)}{(a+1)(kt+i-1)-1}+
\frac{c_1(d_1,d_2-1)}{c_1(d_1,d_2)}\frac{(d_2-1)+k(a-1)}{(a+1)(kt+i-1)-1}+\\+
1-\frac{d_1+d_2+2k(a-1)}{(a+1)(kt+i-1)-1}=
(1-[d_1=k])\frac{d_1-k+ka+a}{(a+1)(kt+i-1)-1}+\\+
(1-[d_2=k])\frac{d_2+k(a-1)}{(a+1)(kt+i-1)-1}+
1-\frac{d_1+d_2+2k(a-1)}{(a+1)(kt+i-1)-1}\le
1+\frac{a}{(a+1)(kt+i-1)-1},
\end{multline*}
we have
\begin{multline}\label{r2inductiont}
\left|\tilde r_2(d_1,d_2,t,i)+\frac{i-1}k([d_2=k]c(d_1)+[d_1=k]c(d_2))t\right|\le\\ \le
C(c_1(d_1,d_2)+c_1(d_2,d_1))\left(kt+\frac{a+\frac12}{a+1}(i-1)\right)\left(1+\frac{a}{(a+1)(kt+i-1)-1}\right)+\\+
[d_1\le k+1]O(c(d_2))+[d_2\le k+1]O(c(d_1)).
\end{multline}
Since
\begin{multline*}
\left(kt+\frac{a+\frac12}{a+1}i\right)-
\left(kt+\frac{a+\frac12}{a+1}(i-1)\right)\left(1+\frac{a}{(a+1)(kt+i-1)-1}\right)=
\frac{a+\frac12}{a+1}-\\-
\left(kt+\frac{a+\frac12}{a+1}(i-1)\right)\frac{a}{(a+1)(kt+i-1)-1}=
\frac1{(a+1)(kt+i-1)-1}\Bigg(\left(a+\frac12\right)(kt+i-1)-\\-
\frac{a+\frac12}{a+1}-akt-a\frac{a+\frac12}{a+1}(i-1)\Bigg)=
\frac1{(a+1)(kt+i-1)-1}\left(\frac{kt}2+\frac{a+\frac12}{a+1}(i-2)\right)
\end{multline*}
is always positive and tends to a nonzero constant limit as $t$ grows,
it is bounded from below by a positive constant. Therefore,
for all sufficiently large values of $C$, the inequality (\ref{r2inductiont})
implies the inductive step by $i$, and so (\ref{r2tildebound}) holds.

As a consequence of (\ref{r2tildebound}), we obtain
$$
\tilde r_2(d_1,d_2,t,i)=O\left(\frac{t}{d_1^{a+1}d_2}+\frac{t}{d_1d_2^{a+1}}\right).
$$

The proven bound, the representation (\ref{covmain}) and Theorem 1 give
the following bound:
$$
\cov(R(d_1,t),R(d_2,t))=O\left(\frac{t}{d_1^{a+1}d_2}+\frac{t}{d_1d_2^{a+1}}\right)
+O(d_1^{-2-a}t)+O(d_2^{-2-a}t)+O\left(\frac1{d_1d_2}\right).
$$
If $d_1\le d_2$, the maximum among the first three terms on the right-hand side is $O(d_1^{-2-a}t)$;
otherwise, the maximum is $O(d_2^{-2-a}t)$. This proves Theorem 2.

\section{Proof of Theorem 3}
We will use the notation $N(s_1,s_2)$ for the number of edges between nodes $s_1$ and $s_2$.
As usual, $N_{t,i}(s_1,s_2)$ is the value of $N(s_1,s_2)$ in the graph before the $i$th step.

First, we define a function
\begin{equation}\label{fdef}
f(d_1,d_2,t,i)=E_{t,i}\left(\sum_{s_1=1}^t\sum_{\substack{s_2=1\\s_2\ne s_1}}^t[\deg s_1=d_1,\deg s_2=d_2]N(s_1,s_2)\right).
\end{equation}
It is easy to see that $EX(d_1,d_2,t)=f(d_1,d_2,t,1)$.

Recurrent equations on $f$ are deduced as it was done in the previous sections. The sum (\ref{fdef})
does not include the last node, so $N(s_1,s_2)$ does not change while adding a new edge.
Thus, the $i$th step acts on $f$ as in the case of $r_2$ (compare with (\ref{r2recur})):
\begin{multline}\label{frecur}
f(d_1,d_2,t,i+1)=f(d_1-1,d_2,t,i)\frac{(d_1-1)+k(a-1)}{(a+1)(kt+i)-1}+\\+
f(d_1,d_2-1,t,i)\frac{(d_2-1)+k(a-1)}{(a+1)(kt+i)-1}+\\+
f(d_1,d_2,t,i)\left(1-\frac{d_1+d_2+2k(a-1)}{(a+1)(kt+i)-1}\right).
\end{multline}

Second, we define a function
\begin{equation}\label{gdef}
g(d_1,d_2,t,i)=E_{t,i}\left([\deg(t+1)=d_2]\sum_{s=1}^t[\deg s=d_1]N(t+1,s)\right).
\end{equation}
Obviously,
\begin{equation}\label{fstart}
f(d_1,d_2,t+1,1)=f(d_1,d_2,t,k+1)+g(d_1,d_2,t,k+1)+g(d_2,d_1,t,k+1)
\end{equation}
and since $N(t+1,s)=0$ before adding any edges from the node $t+1$,
\begin{equation}\label{gstart}
g(d_1,d_2,t,1)=0.
\end{equation}

We now consider one summand of the sum (\ref{gdef}) and the $i$th step.
Let the new edge link nodes $t+1$ and $\gamma$.
We have three non-intersecting cases: $\gamma=s$, $\gamma=t+1$,
$\gamma\not\in\{s,t+1\}$. Note that 
\begin{multline*}
[\gamma=s,\deg_{t,i+1}(t+1)=d_2,\deg_{t,i+1} s=d_1]N_{t,i+1}(s)=\\=
[\gamma=s,\deg_{t,i}(t+1)=d_2-1,\deg_{t,i} s=d_1-1](N_{t,i}(s)+1),
\end{multline*}
\begin{multline*}
[\gamma=t+1,\deg_{t,i+1}(t+1)=d_2,\deg_{t,i+1} s=d_1]N_{t,i+1}(s)=\\=
[\gamma=t+1,\deg_{t,i}(t+1)=d_2-2,\deg_{t,i} s=d_1]N_{t,i}(s),
\end{multline*}
\begin{multline*}
[\gamma\not\in\{s,t+1\},\deg_{t,i+1}(t+1)=d_2,\deg_{t,i+1} s=d_1]N_{t,i+1}(s)=\\=
[\gamma\not\in\{s,t+1\},\deg_{t,i}(t+1)=d_2-1,\deg_{t,i} s=d_1]N_{t,i}(s).
\end{multline*}
Taking the expectation and using the definition (\ref{r2primedef}),
we obtain
\begin{multline}\label{grecur}
g(d_1,d_2,t,i+1)=(g(d_1-1,d_2-1,t,i)+r_2'(d_1-1,d_2-1,t,i))\frac{(d_1-1)+k(a-1)}{(a+1)(kt+i)-1}+\\+
g(d_1,d_2-2,t,i)\frac{(d_2-2)+(i-1)(a-1)+a}{(a+1)(kt+i)-1}+\\+
g(d_1,d_2-1,t,i)\left(1-\frac{(d_1-1)+k(a-1)+(d_2-1)+(i-1)(a-1)+a}{(a+1)(kt+i)-1}\right).
\end{multline}

Third, we derive a bound on $g$. Obviously, $g(d_1,d_2,t,i)=0$ when
$d_2 > 2(i-1)$ or $d_1\ge2(kt+i)$. If $d_1<2(kt+i)$
and $d_2\le 2(i-1)$, then $\frac{d_1+d_2}t\cdot O(d_1^{-a-1})=
O(d_1^{-a}t^{-1})=O(t^{-1})$ and $\frac{d_1+d_2}t\cdot O(t^{-1})=O(t^{-1})$.
Remember that we have proved the bound (\ref{r2prime}) on $r_2'$.
It is easy to see now that
\begin{equation}\label{gasymp}
g(d_1,d_2,t,i+1)=i[d_2=i]c(d_1-1)\frac{(d_1-1)+k(a-1)}{(a+1)k}+O\left(\frac1t\right).
\end{equation}

Finally, we are ready to study $f$.
For the rest of the proof, we will assume that $d_1\ge k$ and $d_2\ge k$.

We denote $A=\frac{\Gamma(a+2)}{\Beta(ka,a+1)}=(a+1)\frac{\Gamma(ka+a+1)}{\Gamma(ka)}$,
$D_1=d_1-k+ka$, $D_2=d_2-k+ka$ for brevity. By definition,
$$c(d_1)=A\frac{\Gamma(D_1)}{\Gamma(D_1+a+2)}.$$

Let $c_X(d_1,d_2)$ be defined recurrently as follows:
\begin{eqnarray*}
c_X(k,k)&=&0,\\
c_X(d_1,k)&=&\frac{(D_1-1)(c_X(d_1-1,k)+c(d_1-1))}{D_1+ka+a+1},\quad d_1>k,\\
c_X(k,d_2)&=&\frac{(D_2-1)(c_X(k,d_2-1)+c(d_2-1))}{D_2+ka+a+1},\quad d_2>k,\\
c_X(d_1,d_2)&=&\frac{(D_1-1)c_X(d_1-1,d_2)+(D_2-1)c_X(d_1,d_2-1)}{D_1+D_2+a+1},\quad d_1,d_2>k.
\end{eqnarray*}

Let
$$
c_2(d_1,d_2)=\frac{\Gamma(D_1)\Gamma(D_2)\Gamma(D_1+D_2+3)}{\Gamma(D_1+2)\Gamma(D_2+2)\Gamma(D_1+D_2+a+2)},
$$
$$
c_3(d_1,d_2)=\frac{\Gamma(D_1)\Gamma(D_2)\Gamma(D_1+D_2+1)}{\Gamma(D_1+1)\Gamma(D_2+1)\Gamma(D_1+D_2+a+2)}.
$$
Obviously, these functions are symmetric. If $d_1>k$,
$$
\frac{c_2(d_1-1,d_2)}{c_2(d_1,d_2)}=\frac{(D_1+1)(D_1+D_2+a+1)}{(D_1-1)(D_1+D_2+2)},
$$
$$
\frac{c_3(d_1-1,d_2)}{c_3(d_1,d_2)}=\frac{D_1(D_1+D_2+a+1)}{(D_1-1)(D_1+D_2)}.
$$
Thus, for $d_1,d_2>k$,
$$
c_2(d_1,d_2)=\frac{(D_1-1)c_2(d_1-1,d_2)+(D_2-1)c_2(d_1,d_2-1)}{D_1+D_2+a+1},
$$
$$
c_3(d_1,d_2)=\frac{(D_1-1)c_3(d_1-1,d_2)+(D_2-1)c_3(d_1,d_2-1)}{D_1+D_2+a+1},
$$

\begin{multline*}
\frac{c(d_1)}{Akac_2(d_1,k)}=\frac{\Gamma(D_1+2)\Gamma(ka+2)\Gamma(D_1+ka+a+2)}{ka\Gamma(D_1+a+2)\Gamma(ka)\Gamma(D_1+ka+3)}=\\=
\frac{ka+1}{D_1+ka+2}\frac{\Gamma(D_1+2)\Gamma(D_1+ka+a+2)}{\Gamma(D_1+a+2)\Gamma(D_1+ka+2)}.
\end{multline*}

Let $(\alpha)_n=\alpha(\alpha+1)\dots(\alpha+n-1)$ be the Pochhammer symbol. Let
${}_2F_1(\alpha,\beta;\gamma;z)=\sum_{n=0}^\infty\frac{(\alpha)_n(\beta)_n}{(\gamma)_n}\frac{z^n}{n!}$
for $\gamma\ne0,-1,-2,\dots$ be the hypergeometric function. According to \cite[15.1.1]{specfunc},
if $\gamma-\alpha-\beta>0$ and $|z|\le1$, this series
converges absolutely. We quote the following formula from \cite[15.1.20]{specfunc}:
$$
{}_2F_1(\alpha,\beta;\gamma;1)=\frac{\Gamma(\gamma)\Gamma(\gamma-\alpha-\beta)}{\Gamma(\gamma-\alpha)\Gamma(\gamma-\beta)}.
$$

Thus
$$
\frac{\Gamma(D_1+2)\Gamma(D_1+ka+a+2)}{\Gamma(D_1+a+2)\Gamma(D_1+ka+2)}=
{}_2F_1(a,ka;D_1+ka+a+2;1)=\sum_{n=0}^\infty\frac{(a)_n(ka)_n}{(D_1+ka+a+2)_n n!}.
$$
Since all terms of the last series are positive and the first term is $1$,
$$
\frac{c(d_1)}{Akac_2(d_1,k)}\ge\frac{ka+1}{D_1+ka+2}.
$$
Moreover,
\begin{multline*}
\frac{\Gamma(D_1+2)\Gamma(D_1+ka+a+2)}{\Gamma(D_1+a+2)\Gamma(D_1+ka+2)}=
1+\sum_{n=0}^\infty\frac{(a)_{n+1}(ka)_{n+1}}{(D_1+ka+a+2)_{n+1}(n+1)!}=\\=
1+\frac{a^2k}{D_1+ka+a+2}\sum_{n=0}^\infty\frac{(a+1)_n(ka+1)_n}{(D_1+ka+a+3)_n(n+1)!} \le \\ \le
1+\frac{a^2k}{D_1+ka+a+2}{\,}_2F_1(a+1,ka+1;2ka+a+3;1) \le \\ \le
1+\frac{a^2k}{D_1+ka+1}\frac{\Gamma(ka+1)\Gamma(2ka+a+3)}{\Gamma(2ka+2)\Gamma(ka+a+2)},
\end{multline*}
$$
\frac{c(d_1)}{Akac_2(d_1,k)}\le\frac{ka+1}{D_1+ka+2}\left(1+\frac{kaB}{D_1+ka+1}\right),
\quad B=a\frac{\Gamma(ka+1)\Gamma(2ka+a+3)}{\Gamma(2ka+2)\Gamma(ka+a+2)}.
$$

In Theorem 3, we have three assertions. The first one says that $ EX(d_1,d_2,t) = c_X(d_1,d_2)t + O_{a,k}(1) $. The second one 
gives a bound for $ c_X $. The third one gives an asymptotic formula for $ c_X $. Now we shall show that our function $ c_X $ 
admits the bound from the second assertion. This bound is equivalent to
\begin{equation}\label{cx}
Aka\left(c_2(d_1,d_2)-\left(4-\frac2{1+ka}\right)c_3(d_1,d_2)\right)\le c_X(d_1,d_2)\le Aka(c_2(d_1,d_2)+Bc_3(d_1,d_2)).
\end{equation}
To prove (\ref{cx}), we use induction on $d_1+d_2$. If $d_1=d_2=k$, the right-hand side of the inequality is obvious,
and its left-hand side follows from
$$
\frac{c_2(k,k)}{c_3(k,k)}=\frac{(2ka+2)(2ka+1)}{(ka+1)(ka+1)}=\frac{4ka+2}{ka+1}=4-\frac{2}{ka+1}.
$$
If $d_1>k$ and $d_2>k$, all the parts of (\ref{cx}) satisfy the same recurrent equation,
so (\ref{cx}) follows from the induction hypothesis. Due to symmetry, it remains
to prove (\ref{cx}) for $d_2=k$, $d_1>k$. We have

\begin{multline*}
\frac{ka\left(c_2(d_1,k)-\frac{2+4ka}{1+ka}c_3(d_1,k)\right)}{kac_2(d_1,k)}=1-\frac{(D_1+1)(2+4ka)}{(D_1+ka+2)(D_1+ka+1)}=\\=
\frac{(D_1-ka)(D_1-ka+1)}{(D_1+ka+1)(D_1+ka+2)}.
\end{multline*}
In particular, $ka\left(c_2(d_1,k)-\frac{2+4ka}{1+ka}c_3(d_1,k)\right)>0$ for $d_1>k$. Then, 
\begin{multline*}
\frac{c_X(d_1,k)}{Aka\left(c_2(d_1,k)-\frac{2+4ka}{1+ka}c_3(d_1,k)\right)}=\\=
\frac{(D_1+ka+1)(D_1+1)(D_1+ka+a+1)}{(D_1-ka)(D_1-ka+1)(D_1-1)}\frac{c_X(d_1,k)}{Akac_2(d_1-1,k)}=\\=
\frac{(D_1+ka+1)(D_1+1)}{(D_1-ka)(D_1-ka+1)}\frac{c_X(d_1-1,k)+c(d_1-1)}{Akac_2(d_1-1,k)} \ge \\ \ge
\frac{(D_1+ka+1)(D_1+1)}{(D_1-ka)(D_1-ka+1)}
\left(\frac{(D_1-1-ka)(D_1-ka)}{(D_1+ka)(D_1+ka+1)}+
\frac{ka+1}{D_1+ka+1}\right)=\\ =
\frac{D_1+1}{D_1-ka+1}\left(\frac{D_1-1-ka}{D_1+ka}+\frac{ka+1}{D_1-ka}\right).
\end{multline*}
Furthermore, 
\begin{multline*}
\frac{D_1-1-ka}{D_1+ka}+\frac{ka+1}{D_1-ka}-\frac{D_1-ka+1}{D_1+1}=
\frac{ka+1}{D_1-ka}-\frac{(ka+1)D_1-k^2a^2+2ka+1}{(D_1+1)(D_1+ka)} \ge \\ \ge
(ka+1)\left(\frac1{D_1-ka}-\frac{D_1+ka+1}{(D_1+1)(D_1+ka)}\right)=
(ka+1)ka\frac{D_1+ka+2}{(D_1+1)(D_1+ka)(D_1-ka)}>0.
\end{multline*}
This proves the left hand-side of the inequality (\ref{cx}).

Now, 
\begin{multline*}
\frac{c_X(d_1,k)}{Aka(c_2(d_1,k)+Bc_3(d_1,k))}=\\=
\left(1+B\frac{(D_1+1)(ka+1)}{(D_1+ka+1)(D_1+ka+2)}\right)^{-1}
\frac{(D_1+1)(D_1+ka+a+1)}{(D_1-1)(D_1+ka+2)}
\frac{c_X(d_1,k)}{Akac_2(d_1-1,k)} \le \\ \le
\frac{(D_1+ka+1)(D_1+1)}{(D_1+ka+1)(D_1+ka+2)+B(D_1+1)(ka+1)}\times \\ \times
\left(1+B\frac{D_1(ka+1)}{(D_1+ka)(D_1+ka+1)}+\frac{ka+1}{D_1+ka+1}\left(1+\frac{kaB}{D_1+ka}\right)\right).
\end{multline*}
Furthermore, 
\begin{multline*}
\left(1+B\frac{D_1(ka+1)}{(D_1+ka)(D_1+ka+1)}+\frac{ka+1}{D_1+ka+1}\left(1+\frac{kaB}{D_1+ka}\right)\right)-\\-
\left(\frac{D_1+ka+2}{D_1+1}+B\frac{ka+1}{D_1+ka+1}\right)=
\left(\frac{ka+1}{D_1+ka+1}-\frac{ka+1}{D_1+1}\right)<0.
\end{multline*}
This proves the right-hand side of the inequality (\ref{cx}).

The third assertion of Theorem 3 (i.e., the asymptotic formula for $c_X(d_1,d_2)$) 
is derived from (\ref{cx}) in the same way as a similar formula was derived from (\ref{casymp}).

It remains to prove the first assertion. We will use the bound $c_X(d_1,d_2)=O\left(\frac{(d_1+d_2)^{1-a}}{d_1^2d_2^2}\right)$.

Let
$$
\tilde f(d_1,d_2,t,i)=f(d_1,d_2,t,i)-c_X(d_1,d_2)\left(t+\frac ik-\frac1{k(a+1)}\right).
$$

Let $C=C(a,k)$ be a sufficiently large constant which will be determined later.
We claim that
\begin{multline}\label{ftildebound}
\left|\tilde f(d_1,d_2,t,i)+(i-1)\left([d_1=k]\frac{(D_2-1)c(d_2-1)}{(a+1)k}+[d_2=k]\frac{(D_1-1)c(d_1-1)}{(a+1)k}\right)\right|\le\\
\le C\left(1-\frac1{(a+2)k(t+1)}\right)^{i-1}.
\end{multline}

Since $r_2'(d_1,d_2,t,i)=0$ when $d_1+d_2>kt+k+2(i-1)$, (\ref{gstart}) and
(\ref{grecur}) imply that $g(d_1,d_2,t,i)=0$ when $d_1+d_2>kt+k+2(i-1)$.
Consequently, (\ref{frecur}) and (\ref{fstart}) imply that
$f(d_1,d_2,t,i)=0$ when $d_1+d_2>kt+2k+(i-1)$.

If $d_1+d_2>kt+2k+(i-1)$, then
$$\tilde f(d_1,d_2,t,i)=-c_X(d_1,d_2)\left(t+\frac ik-\frac1{k(a+1)}\right)=O\left(\frac{(d_1+d_2)^{2-a}}{d_1^2d_2^2}\right),$$
so (\ref{ftildebound}) holds for all sufficiently large values of $C$.

Now assume $d_1+d_2\le kt+2k+(i-1)$. We will prove (\ref{ftildebound}) by induction
on $t$ and, for fixed $t$, on $i$. The basis of induction $t=1,\dots,2+\lfloor\frac1{ka}\rfloor$ and any $i=1,\dots,k+1$
obviously holds for all sufficiently large values of $C$.

Now let $t\ge 3+\lfloor\frac1{ka}\rfloor$ and let (\ref{ftildebound}) hold for $t-1$. Since (\ref{ftildebound})
is trivial for $d_1=d_2=k$ and symmetrical, we may assume that $d_1>k$.
From (\ref{fstart}) and (\ref{gasymp}), we obtain
\begin{multline*}
\left|\tilde f(d_1,d_2,t,1)\right|=\left|\tilde f(d_1,d_2,t-1,k+1)+k[d_2=k]c(d_1-1)\frac{D_1-1}{(a+1)k}+O\left(\frac1t\right)\right|\\
\le C\left(1-\frac{1}{(a+2)kt}\right)^k+O\left(\frac1t\right).
\end{multline*}
The right-hand side is less than $C$ for all sufficiently large values of $C$.

Finally, let $t\ge3+\lfloor\frac1{ka}\rfloor>2+\frac1{ka}$, $i>1$ and let (\ref{ftildebound}) hold for $i-1$. We reuse the notation
$T=(a+1)(kt+i-1)-1$ and again assume $d_1>k$. From (\ref{frecur}) we obtain
\begin{multline*}
\tilde f(d_1-1,d_2,t,i-1)\frac{D_1-1}{T}+[d_2>k]\tilde f(d_1,d_2-1,t,i-1)\frac{D_2-1}T+\\+
\tilde f(d_1,d_2,t,i-1)\left(1-\frac{D_1+D_2}T\right)=
\tilde f(d_1,d_2,t,i)+c_X(d_1,d_2)\left(\frac{T}{(a+1)k}+\frac1k\right)-\\-
\frac{c_X(d_1-1,d_2)(D_1-1)+[d_2>k]c_X(d_1,d_2-1)(D_2-1)+c_X(d_1,d_2)(T-(D_1+D_2))}{(a+1)k}=\\=
\tilde f(d_1,d_2,t,i)+[d_2=k]c(d_1-1)\frac{D_1-1}{(a+1)k}.
\end{multline*}
The assumptions $d_1+d_2\le kt+2k+(i-1)$ and $t>2+\frac1{ka}$ imply that
$1-\frac{D_1+D_2}{T}=1-\frac{d_1+d_2+2k(a-1)}{(a+1)(kt+i-1)-1}\ge0$.
Since $(D_1-2)c(d_1-2)=[d_1>k+1](D_1+a)c(d_1-1)$,
\begin{multline*}
(i-2)[d_2=k]\frac{(D_1-2)c(d_1-2)}{(a+1)k}\frac{D_1-1}T+(i-2)[d_2=k]\frac{(D_1-1)c(d_1-1)}{(a+1)k}\left(1-\frac{D_1+D_2}T\right)=\\=
(i-2)[d_2=k]\frac{(D_1-1)c(d_1-1)}{(a+1)k}\left([d_1>k+1]\frac{D_1+a}T+1-\frac{D_1+D_2}T\right)=\\=
(i-2)[d_2=k]\frac{(D_1-1)c(d_1-1)}{(a+1)k}+O\left(\frac1t\right).
\end{multline*}
If $d_2>k$, then $\frac{D_1-1}T+\frac{D_2-1}T+1-\frac{D_1+D_2}T=1-\frac2T\le1-\frac1T$. If
$d_2=k$, then $\frac{D_1-1}T+1-\frac{D_1+D_2}T=1-\frac{ka+1}T\le1-\frac1T$.
Thus,
$$
\left|\tilde f(d_1,d_2,t,i)+(i-1)[d_2=k]\frac{(D_1-1)c(d_1-1)}{(a+1)k}\right|\le
$$
$$
\le
C\left(1-\frac1{(a+2)k(t+1)}\right)^{i-2}\left(1-\frac1{(a+1)k(t+1)}\right)+O\left(\frac1t\right).
$$
For all sufficiently large values of $C$, the right-hand side is less than $C\left(1-\frac1{(a+2)k(t+1)}\right)^{i-1}$,
so the induction on $i$ is complete.

Theorem 3 follows from the proven bound (\ref{ftildebound}).

\section{Proof of Theorem 4}
We use the Azuma--Hoeffding inequality.
\begin{theorem} \cite{azuma}, \cite{hoeffding} Let $(X_s)_{s=0}^n$ be a martingale with
$|X_{s+1}-X_s|\le\delta$ for $s=0,\dots,n-1$, and $x>0$. Then
$$
P\left(|X_n-X_0|\ge x\right)\le2\exp\left(-\frac{x^2}{2c^2n}\right).
$$
\end{theorem}

We fix $d_1,d_2,t$ and denote $X=X(d_1,d_2,t)$. Let $G$ be a random graph in $H_{a,k}^{(t)}$; it has $kt$ edges,
sorted by the creation time. Let $G^{(s)}$ be a graph with $s$ first edges.
Let $X_s=E(X|G^{(s)}), s=0,\dots,kt$. In this sequence $X_0=EX$, $X_{kt}=X$.
By definition of the probabilistic space, the sequence $X_s$ is a martingale.
We will estimate possible differences between adjacent elements of the sequence.

We fix any $s$ from $0$ to $kt-1$. Let $v$ be the head of the last edge in $G^{(s+1)}$,
so $v$ is a random quantity depending on $G$. By definition
$$\begin{array}{rcrl}
X_s &=& \sum_\gamma \Pr(v=\gamma) & E(X|G^{(s)},v=\gamma),\\
X_{s+1} &=& & E(X|G^{(s)},v=v(G^{(s+1)})),
\end{array}$$
where the sum is over all nodes of $G$. Hence it is clear that
$$
\min_{\gamma} E(X|G^{(s)},v=\gamma)\le X_s,X_{s+1}\le\max_\gamma E(X|G^{(s)},v=\gamma),
$$
$$
|X_s-X_{s+1}|\le\max_\gamma E(X|G^{(s)},v=\gamma)-\min_\gamma E(X|G^{(s)},v=\gamma).
$$
Let $\gamma_1\in\arg\min E(X|G^{(s)},v=\gamma)$ and $\gamma_2\in\arg\max E(X|G^{(s)},v=\gamma)$.
It is sufficient to prove an upper bound for
$$
E(X|G^{(s)},v=\gamma_2)-E(X|G^{(s)},v=\gamma_1).
$$
We consider the sum
\begin{equation}\label{xsum}
X=\sum_{s_1=1}^t\sum_{\substack{s_2=1\\s_2\ne s_1}}^t[\deg s_1=d_1,\deg s_2=d_2]N(s_1,s_2).
\end{equation}
Replacing the condition
$v=\gamma_1$ by the condition $v=\gamma_2$ changes
distributions of degrees of $\gamma_i$ and distributions
of $N(\gamma_i,*)=N(*,\gamma_i)$; distributions of other
values of $N$ do not change. Thus distributions of all
terms in the sum (\ref{xsum}) except those with
$\{\gamma_1,\gamma_2\}\cap\{s_1,s_2\}\ne\varnothing$ are the same for
$v=\gamma_1$ and $v=\gamma_2$. Let
$$
X'=\sum_{s_1=1}^t\sum_{\substack{s_2=1\\s_2\ne s_1\\ \{s_1,s_2\}\cap\{\gamma_1,\gamma_2\}\ne\varnothing}}^t
[\deg s_1=d_1,\deg s_2=d_2]N(s_1,s_2).
$$
Then
$$
E(X-X'|G^{(s)},v=\gamma_1)=E(X-X'|G^{(s)},v=\gamma_2).
$$
Obviously, $X'\ge0$. We have
\begin{multline*}
X'\le\sum_{s_1=1}^t[\deg s_1=d_1,\deg\gamma_1=d_2]N(s_1,\gamma_1)+\sum_{s_1=1}^t[\deg s_1=d_1,\deg\gamma_2=d_2]N(s_1,\gamma_2)+\\+
\sum_{s_2=1}^t[\deg\gamma_1=d_1,\deg s_2=d_2]N(\gamma_1,s_2)+\sum_{s_2=1}^t[\deg\gamma_2=d_1,\deg s_2=d_2]N(\gamma_2,s_2)\le\\ \le
[\deg\gamma_1=d_2]\sum_{s_1=1}^t N(s_1,\gamma_1)+[\deg\gamma_2=d_2]\sum_{s_1=1}^t N(s_1,\gamma_2)+\\+
[\deg\gamma_1=d_1]\sum_{s_2=1}^t N(\gamma_1,s_2)+[\deg\gamma_2=d_1]\sum_{s_2=1}^t N(\gamma_2,s_2)=\\=
[\deg\gamma_1=d_2]d_2+[\deg\gamma_2=d_2]d_2+[\deg\gamma_1=d_1]d_1+[\deg\gamma_2=d_1]d_1\le2(d_1+d_2).
\end{multline*}
Thus,
$$
0\le E(X'|G^{(s)},v=\gamma_1),E(X'|G^{(s)},v=\gamma_2)\le2(d_1+d_2),
$$
$$
|E(X'|G^{(s)},v=\gamma_1)-E(X'|G^{(s)},v=\gamma_2)|\le2(d_1+d_2),
$$
$$
|X_s-X_{s+1}|\le2(d_1+d_2).
$$
Consequently, the sequence $(X_s)$ satisfies the condition of Theorem 5 with $n=kt$ and $\delta=2(d_1+d_2)$.
Substituting $x=c(d_1+d_2)\sqrt{kt}$ in Theorem 5, we obtain Theorem 4.

\end{document}